\documentclass[letterpaper, 10 pt, conference]{ieeeconf}

\usepackage{amssymb}
\usepackage{amsfonts}
\usepackage{cite}
\usepackage{graphicx}
\usepackage{psfrag}
\usepackage{subfigure}
\usepackage{stfloats}
\usepackage{amsmath}
\usepackage{bbm}
\usepackage{array}
\IEEEoverridecommandlockouts

\def\Rset{\mathbb{R}}

\def\diag{\operatorname{diag}}

\title{\LARGE \bf
Pattern Recognition Facilities of Extended Kalman Filtering in Stochastic Neural Fields}

\author{Maria V.~Kulikova$^{1}$, Pedro M.~Lima$^{1}$, Gennady Yu.~Kulikov$^{1}$
\thanks{*This work was supported by the Portuguese FCT~--- \emph{Funda\c{c}\~ao para a Ci\^encia e a Tecnologia}, via the projects UIDB/04621/2020, UIDP/04621/2020, and PTDC/MAT-APL/31393/2017.}
\thanks{$^{1}$The authors are with CEMAT (Center for Computational and Stochastic Mathematics), Instituto Superior T\'{e}cnico, Universidade de Lisboa,
          Av. Rovisco Pais 1,  1049-001 LISBOA, Portugal;
        {\tt maria.kulikova at ist.utl.pt; pedro.t.lima at ist.utl.pt; gennady.kulikov at tecnico.ulisboa.pt}}%
}

\begin{document}

\maketitle
\thispagestyle{empty}
\pagestyle{empty}

\begin{abstract}

In mathematical neuroscience, a special interest is paid to a working memory mechanism 
in the neural tissue modeled by the Dynamic Neural Field (DNF) in the presence of model uncertainties. 
The working memory facility implies that the neurons' activity remains self-sustained after the external 
stimulus removal due to the recurrent interactions in the networks and allows the system to cope with 
missing sensors' information. In our previous works, we have developed two reconstruction methods of 
the neural membrane potential from {\it incomplete} data available from the sensors. The methods are 
derived within the Extended Kalman filtering approach by using the Euler-Maruyama method and the It\^{o}-Taylor expansion of order 1.5. It was shown that the It\^{o}-Taylor EKF-based restoration process is more accurate than the Euler-Maruyama-based alternative. It improves the membrane potential reconstruction quality in case of incomplete sensors information. The aim of this paper is to investigate their pattern recognition facilities, i.e. the quality of the pattern formation reconstruction in case of model uncertainties and incomplete information. The numerical experiments are provided for an example of the stochastic DNF with multiple active zones arisen in a neural tissue.
\end{abstract}

\section{INTRODUCTION}\label{sect1}

In~\cite{Am77}, the conditions for the existence of locally excited regions, called bumps, in the Dynamic Neural Fields (DNFs) have been explored. The DNF equation suggested in the cited paper is a mathematical description of cortical neural tissue modeled by a partial integro-differential equation. In mathematical neuroscience, a special interest is paid to the scenario when the activity remains self-sustained after the external stimulus removal due to the recurrent interactions in the networks, i.e. to the stability of the bump regions. This allows the working memory mechanism to be investigated and implemented in robotics. The  working memory facility is defined as the capacity to transiently hold sensory information to guide forthcoming actions~\cite{lima2018numerical}. The mathematical analysis in~\cite{Am77} concerns with the existence of one-bump solution, which is then extended in~\cite{FeEr16} to the derivation of conditions of the existence and stability of multi-bump solutions representing the memory of a series of transient inputs.

In the past few years, there is a growing body of literature that recognizes the importance of exploring the DNFs in the presence of model uncertainties. This yields the Stochastic Dynamic Neural Field (SDNF) models given as follows~\cite{KuRi14}:
\begin{align}
du(x,t) & =\biggl[I(x,t) + \int_{\Omega} F(|x-y|) S(u(y,t)-\theta)dy    \nonumber \\
& - \alpha u(x,t) \biggr]dt + \epsilon dW(x,t), \; t\in [0, T] \label{eq1.3}
\end{align}
where $\alpha\in\Rset$, $\alpha > 0$ is the decay rate, $u(x,t)$ is the average membrane potential at time $t$ of a neuron located at position $x \in\Omega  \subset \Rset$. Its average activity (i.e. the firing rate) is given by the transfer function $S(\cdot)$, which is often assumed to be the Heaviside function. More precisely, the Heaviside function $H(u(x,t)-\theta)$ is defined by $H(\cdot)=1$ when $u(x,t) \ge \theta$ and $H(\cdot)=0$ when $u(x,t) < \theta$ at time $t$. The parameter $\theta$ is the default level to which the potential $u(x,t)$ is relaxed when the external inputs are absent. In other words, the firing rate function with threshold $\theta \ge 0$ maps the synaptic signal, $u(x,t)$, to a resulting fraction of active neurons. Alternatively, it might be considered as a probability of activation of a single neuron. Furthermore, each neuron $x$ is affected by the time-varying external stimulus $I(x,t)$, i.e. the time-varying input from the outside of the field.  Besides, there is an impact through lateral connections from each of its field neighbors. This effect is modeled by the connectivity function $F(\cdot)$. As can be seen, the average strength of the connection between any two neurons $x$ and $y$ is assumed to be a function of their distance, only. This means that the interaction is symmetric and the neural field is homogeneous. Additionally, we explore the symmetric spatial domain with the length $2L$, i.e. $x\in\Omega\equiv[-L,L]$. In addition, the scalar $\epsilon\ge 0$ stands for the noise level in our model and $W(\cdot)$ is a trace-class space-valued $Q$-Wiener process. The partial integro-differential equation in~\eqref{eq1.3} should be solved with the given initial value $u(x,0)=u_0(x)$ where the function $u(x,t)$ is assumed to satisfy the periodic boundary conditions on the spatial domain $\Omega$.

Following~\cite{KuRi14}, if the function $S(\cdot)$ is assumed to be globally Lipschitz continuous on $\Rset$ and $F(\cdot)$ is such that the integral operator in equation~\eqref{eq1.3} is self-adjoint and compact, then it ensures the existence and uniqueness of a mild solution to the SDNF equation in~\eqref{eq1.3} for any sufficiently smooth function $u_0(x)$.  The solution $u(x,t)$ of equation~\eqref{eq1.3} is expected to have a symmetric form in case of the homogeneous field examined, the symmetric domain $\Omega$,  the periodic boundary conditions on $\Omega$ and the symmetric external stimulus $I(x,t)$. This are the common assumptions in mathematical neuroscience literature; e.g., see~\cite{FeEr16,LiBu15}.

A considerable amount of literature has been published on numerical simulation and theoretical research devoted to the conditions of the existence and stability of one- and  multi-bump solutions in the SDNF equations given by~\eqref{eq1.3}. In this paper, we explore the reconstruction problem of the SDNF solution from the incomplete information available from the sensors. In particular, we are interested in the recognition of the pattern formation in a neural tissue. The membrane potential is assumed to be partially observable in both domains, which are the time and spatial domains. It is widely acknowledged that the self-sustained properties of the neural population dynamics allow the system to cope with missing sensory information and to anticipate the action outcomes ahead of their realization. The goal of this paper is to explore the pattern recognition facilities of the reconstruction methods developed in~\cite{2021:Kulikova:Romania,2022:Kulikova:ECC}, i.e. their restoration quality of the pattern formation in a neural tissue in case of model uncertainties and incomplete data.

\section{BRIEF DESCRIPTION OF STATE-SPACE REPRESENTATION AND METHODS}\label{sect1}

To solve the reconstruction problem by the filtering methods in~\cite{2021:Kulikova:Romania,2022:Kulikova:ECC}, we discuss how the state-space model can be obtained for the SDNF under examination. Following~\cite{KuRi14,KuLi21CAM}, we utilize the Karhunen-Lo\`eve formula for approximating the solution $u(x,t)$  as follows~\cite{Lo78}:
\begin{equation} \label{eq2.8}
u(x,t)=\sum_{k=0}^{K} u_k(t) v_k(x)
\end{equation}
where $u_k(t)$ are pairwise uncorrelated random processes defined on the time interval $[0,T]$ and $v_k(x)$ are deterministic continuous real-valued functions defined on the spatial domain $\Omega$. The functions $v_k(x)$ are the eigenfunctions of the covariance operator in the stochastic noise term of the SDNF equation in~\eqref{eq1.3}. Following~\cite[Proposition~2.1.10]{PrRo07}, the eigenfunctions $v_k(x)$ are utilized in the following expansion:
\begin{equation}\label{eq2.4}
d W(x,t)=\sum_{l=0}^{K} v_l(x) \lambda_l  d\beta_l(t)
\end{equation}
where the stochastic processes $\beta_l(t)$ stand for mutually independent Brownian motions with zero mean and unit variance, and
 $\lambda_l$ are the eigenvalues of the covariance operator.

In this paper, we explore the scenario discussed in~\cite[Example~2.1]{KuRi14} where the eigenvalues, $\lambda_l$, $l = 0, \ldots, K$, of the covariance operator are computed as follows~\cite{Sh05}:
\begin{equation}
\lambda_l = \sqrt{\exp\bigl\{-\xi^2 l^2/(4 \pi)\bigr\}} \label{eq:lambdas}
\end{equation}
where $\xi$ is fixed and denotes the spatial correlation length.

Next, the trigonometric orthonormal basis is utilized (the sinusoidal eigenfunctions are absent because of the even form of the SDNF solution under examination), i.e. we have
\begin{equation} \label{eq2.9}
v_k(x) = \left\{
\begin{array}{ll}  1/\sqrt{2L}, & \mbox{if}\quad k=0,\\
             1/\sqrt{L}\cos\{k \pi x/L\}, & \mbox{if}\quad k\ge 1.
            \end{array}
\right.
\end{equation}

The equidistant mesh is introduced on the spatial domain $\Omega = [-L,L]$ with the step size $h_x=2L/N$, i.e.
 \[\mathbf{x}:=\left\{x_0=-L,\;x_{i}=x_{0}+i h_x,\; i=0,1,\ldots,N,\;x_N=L\right\}\]
where $N$ is some fixed number of subdivisions (the condition $N>>K$ should be satisfied to get an accurate approximation). Each function $v_k(x)$, $k=0, \ldots, K$ in~\eqref{eq2.9} is calculated on the mesh $\mathbf{x} \in {\mathbb R}^{N+1}$ and denoted by $\mathbf{v}_k:=v_k(\mathbf{x})$.

Following~\cite{KuLi21CAM}, we utilize the Trapezoidal Rule for calculating the integrals in equation~\eqref{eq1.3} and derive a set of stochastic differential equations (SDEs) as follows:
\begin{align}
& d\mathbf{u}(t) =  \biggl[h_x\,\mathbf{V}\times I(\mathbf{x},t) - \alpha \mathbf{u}(t) +\mathbf{FS}(u(\mathbf{x},t),\theta)\biggr]dt \nonumber \\
 &   + \epsilon \mathbf{\Lambda} d\mathbf{\beta}(t), \; \mathbf{u}(t):=[u_0(t),u_1(t),\ldots,u_K(t)]^\top \label{eq2.13}
\end{align}
where $\beta(t)$ is the $(K+1)$-dimensional standard Brownian motion whose increment $d\mathbf{\beta}(t)$ is independent of $\mathbf{u}(t)$ and $\mathbf{\Lambda}:=\diag\{\lambda_0,\lambda_1,\ldots,\lambda_K\}$ with each $\lambda_k$ computed by~\eqref{eq:lambdas}. We examine the Gaussian noise scenario, i.e. the increment $d\mathbf{\beta}(t)$ is a standard Gaussian white process with zero-mean and the covariance $Q\,dt>0$ where $Q=I_{K+1}$ and $I_{K+1}$ is an identity matrix of size $K+1$. The readers are referred to~\cite{KuLi21CAM} for more details and derivation. Here, we stress that the matrix $\mathbf{V}_f$ is formed from the vectors $\mathbf{v}_k:=v_k(\mathbf{x})$ by rows and, hence, it has the dimension of $(K+1)\times (N+1)$. Meanwhile, due to the even form of the SDNF solution and the Trapezoidal Rule utilized for computing the integrals, we also have
$\mathbf{V} \in {\mathbb R}^{(K+1)\times N}$, i.e. the matrix $\mathbf{V}_f$ without the last column. It is utilized for calculating the following term:
 \begin{align}
\mathbf{FS}(u(\mathbf{x},t),\theta) & = h_x^2\, \mathbf{V} \times \mathbf{F}\times \mathbf{s} \label{eq2.14}
\end{align}
where the symmetric matrix $\mathbf{F}$ consists of the entries $F(|x_j-x_i|)$, $i,j=0,1,\ldots,N-1$. The column-vector $\mathbf{s}$ of size $N$ contains the entries, which are calculated by the formula $\mathbf{s}_i = S(\mathbf{u}^\top(t)\mathbf{V}_i - \theta)$ where $\mathbf{V}_i$ means the $i$th column in the matrix $\mathbf{V}$, which corresponds to all basis functions $v_k(x)$, $k=0,1,\ldots,K$, evaluated at the spatial node $x_i$, $i=0,\ldots, N-1$. Finally, the $N$-dimensional vector $I(\mathbf{x},t)$ in equation~\eqref{eq2.13} is the external input function calculated on the mesh $\mathbf{x}$ with entries $x_i$ and $i = 0,1,\ldots,N-1$.

Equation~\eqref{eq2.13} represents the process equation of the standard continuous-discrete nonlinear state-space model, i.e.
\begin{align}
d\mathbf{u}(t) & = f\bigl(t,\mathbf{u}(t)\bigr)dt + G d\mathbf{\beta}(t), \quad t>0,  \label{eq1.1:new}
\end{align}
where $\mathbf{u}(t)$ is the unknown state vector to be estimated and the diffusion matrix equals to $G:=\epsilon \mathbf{\Lambda}$. The vector-function $f:{\mathbb R}\times{\mathbb R}^{(K+1)}\to {\mathbb R}^{(K+1)}$ is the drift function given by
\begin{align}
f(\cdot) & = \biggl[h_x\,\mathbf{V}\times I(\mathbf{x},t) - \alpha \mathbf{u}(t) +\mathbf{FS}(u(\mathbf{x},t),\theta)\biggr]. \label{drift:state:fun}
\end{align}

The measurement equation of the standard state-space model can be found from the Karhunen-Lo\`eve formula~\eqref{eq2.8}. Indeed, the solution $u(x,t)$ of the SDNF equation~\eqref{eq1.3} is reconstructed from the data partially observed over the time domain $[0, T]$ with the sampling rate (sampling period) $\Delta_t=t_{k}-t_{k-1}$ through the measurement equation as follows:
\begin{align}
z(t_k):=u({\mathbf x},t_k)  & =  \mathbf{V}_{f}^{\top}\mathbf{u}(t_k)+\eta(t_k), k =1,2,\ldots \label{eq1.2:new}
\end{align}
where the matrix $\mathbf{V}_{f}$ is of dimension $(K+1)\times (N+1)$ as explained in detail after formula~\eqref{eq2.14}. More precisely, the $(k,i)$-entry of the matrix  $\mathbf{V}_{f}$ means the value of the $k$-th eigenfunction~\eqref{eq2.9} computed at the space node $x_i$, i.e. $v_k(x_i)$, $k=0,1,\ldots,K$, $i=0,1,\ldots,N$. We also assume that $\eta(t_k)$ is a zero-mean Gaussian noise with the covariance $R>0$.

In this paper, we examine the reconstruction problem when the incomplete data set additionally appears over the spatial domain $\Omega = [-L,L]$. This means that a few locations ${\mathbf x}_{\{i\}}$ are supplied by the sensors for collecting the information about the membrane potential, only. The set ${\mathbf x}_{\{i\}}$ denotes a subset of the mesh ${\mathbf x}\in {\mathbb R}^{N+1}$ where each position $x_{\{i\}}$ is supplied by the sensor, which collects the data, i.e. the measurements are taken at ${\mathbf z}_{\{i\}}(t_k)$. Thus, each row of the matrix $\mathbf{V}_{f}$ has missing entries. Having denoted the number of sensors by $m$, $m \le N+1$, we get
\begin{equation}
{\mathbf z}_{\{i\}}(t_k)  = [\mathbf{V}_{f}^{\top}]_{\{i\}}\mathbf{u}(t_k) +\eta(t_k) \label{eq:meas}
\end{equation}

To summarize, the nonlinear filtering problem associated with the SDNF-oriented state-space model in equations~\eqref{eq2.13}, \eqref{eq:meas} permits the reconstruction process of the membrane potential given the incomplete information over $[0, T]\times [-L, L]$. The filter should utilize the initial state value that is also assumed to be normally distributed with the mean ${\mathbf u}(t_0)$ and covariance $\Pi_0>0$. The average membrane potential ${\mathbf u}(t_0)$ is defined from the  Karhunen-Lo\`eve formula~\eqref{eq2.8} for approximation of the initial value $u(x,0)=u_0(x)$ of the SNFE in~\eqref{eq1.3}, i.e. $u(\mathbf{x},t)=V_f^{\top}\mathbf{u}(t)$. Thus, given the initial value $u_0(\mathbf{x},t)$ of the SNFE discretized on the mesh $\mathbf{x}$, we may find the initial state $\mathbf{u}(t_0)$ as follows:
\begin{equation} \label{SDE:initial}
\mathbf{u}(t_0) = (V_f \; V_f^{\top})^{-1} V_f u_0(\mathbf{x},t).
\end{equation}

To solve the stated filtering problem, we suggest to utilize the extended Kalman filter (EKF) because the measurement equation~\eqref{eq:meas} in the state-space model has a linear form.

\subsection{The Euler-Maruyama-based EKF reconstruction method}

The Euler-Maruyama method for solving the SDNF in~\eqref{eq1.3} has been designed in~\cite{KuLi21CAM}, meanwhile the related EKF estimator for reconstructing the membrane potential through the state-space representation has been proposed in~\cite{KuKu19cIEEE_ICSTCC}. We briefly explain the main steps.

\textsc{Time Update}: On each sampling interval $\Delta_t=t_{k}-t_{k-1}$, introduce the mesh $t_{k-1}^{(l)} = t_{k-1}\!+l\delta$, $l=0,\ldots,L-1$. The filtered estimate $\mathbf{\hat u}_{k-1|k-1}$ and the error covariance $P_{k-1|k-1}$ are available from the previous recursion step. Set the initial values for the current prediction step: $\mathbf{\hat u}^{(0)}_{k|k-1}:=\mathbf{\hat u}_{k-1|k-1}$ and $P^{(0)}_{k|k-1}:=P_{k-1|k-1}$. For $l=0,\ldots,L-1$, compute
\begin{align}
\mathbf{\hat u}_{k|k-1}^{(l+1)}  & = \mathbf{\hat u}_{k|k-1}^{(l)}  + \delta f\bigl(t_{k-1}^{(l)},\mathbf{\hat u}_{k|k-1}^{(l)}\bigr) \label{EM:EKF:mean} \\
P^{(l+1)}_{k|k-1} & = [I - \delta J\bigl(t_{k-1}^{(l)},\mathbf{\hat u}_{k|k-1}^{(l)}\bigr)] \nonumber \\
& \times P^{(l)}_{k|k-1}[I - \delta J\bigl(t_{k-1}^{(l)},\mathbf{\hat u}_{k|k-1}^{(l)}\bigr)]^\top \!\!+ \delta GG^\top \label{EM:EKF:cov}
\end{align}
where the Jacobian matrix $J = \partial f\bigl(t,\mathbf{u}\bigr)/\partial \mathbf{u}$ can be found from equation~\eqref{drift:state:fun} as follows:
\begin{align}
J = -\alpha I_{K+1} + h_x^2\, \mathbf{V} \times \mathbf{F}\times \partial \mathbf{s}/\partial \mathbf{u} \label{Jacobian:drift}
\end{align}
where each $i$th ($i = 1, \ldots,K+1$) row of the Jacobian matrix $\partial \mathbf{s}/\partial \mathbf{u}$ is computed in the following way:
\begin{align}
\frac{\partial}{\partial \mathbf{u}}\mathbf{s}_i & = S^{\prime}\bigl(\mathbf{u}^\top\mathbf{V}_i - \theta\bigr) \times \mathbf{V}^\top_i
\end{align}
where the prime refers to the derivative of the firing rate function $S(\cdot)$ evaluated at the argument $\mathbf{u}^\top\mathbf{V}_i - \theta$ and, as customary, the vector $\mathbf{V}_i$ is the $i$th column in the matrix $\mathbf{V}$.

\textsc{Measurement Update}: Define $H:=[\mathbf{V}_{f}^{\top}]_{\{i\}}$ and update the estimate as follows:
\begin{align}
 R_{e,k}  & = HP_{k|k-1}H^{\top} + R, \quad K_k =  P_{k|k-1}H^{\top} R_{e,k}^{-1}, \label{kf:f:K} \\
\mathbf{\hat u}_{k|k} & =  \mathbf{\hat u}_{k|k-1} + K_k \left[{\mathbf z}_{\{i\}}(t_k) - H \mathbf{\hat u}_{k|k-1}\right], \label{kf:f:ek} \\
P_{k|k} & =  (I -K_kH)P_{k|k-1}. \label{kf:f:P}
\end{align}

Having computed the estimates $\mathbf{\hat u}_{k|k}$, the membrane potential, $u(x,t)$, which obeys SDNF equation~\eqref{eq1.3}, can be reconstructed from the incomplete information at each time instance $t_k$ as follows: $\mathbf{u}(\mathbf{x}_{\{i\}},t_k) = [\mathbf{V}_{f}^{\top}]_{\{i\}}\mathbf{\hat u}_{k|k}$. A fast reconstruction procedure based on the {\it sequential} measurement update has been proposed in~\cite{kulikova2022sequential}.

\subsection{The It\^{o}-Taylor-based EKF reconstruction method}

Recall, the Euler-Maruyama method is of strong order 0.5 (EM-0.5). The EKF estimator can be constructed within a higher order method, for instance, by using the It\^{o}-Taylor expansion of order 1.5 (IT-1.5); see~\cite{1999:Kloeden:book}. The IT-1.5 integrator for solving the SDNF  equation has been developed in~\cite{KuKu20IEEE_ICSTCC}, meanwhile the related EKF estimator for reconstructing the membrane potential has been recently proposed in~\cite{2022:Kulikova:ECC}. Here, we briefly discuss the main steps of the estimation method.

\textsc{Time Update}: On each sampling interval $\Delta_t=t_{k}-t_{k-1}$, introduce the mesh $t_{k-1}^{(l)} = t_{k-1}\!+l\delta$, $l=0,\ldots,L-1$. The filtered estimate $\mathbf{\hat u}_{k-1|k-1}$ and the error covariance $P_{k-1|k-1}$ are available from the previous recursion step. Set the initial values for the current prediction step: $\mathbf{\hat u}^{(0)}_{k|k-1}:=\mathbf{\hat u}_{k-1|k-1}$ and $P^{(0)}_{k|k-1}:=P_{k-1|k-1}$. For $l=0,\ldots,L-1$, compute
\begin{equation}
\mathbf{\hat u}_{k|k-1}^{(l+1)} =  f_d\bigl(t_{k-1}^{(l)},\mathbf{\hat u}_{k|k-1}^{(l)}\bigr) \label{EKF:mean}
\end{equation}
where the discretized drift function is given as follows~\cite{1999:Kloeden:book}:
\begin{align}
\!\!\!\!f_d(\cdot)\! =\! \mathbf{u}_{k-1}^{(l)}\!+\delta f\bigl(t_{k-1}^{(l)},\mathbf{u}_{k-1}^{(l)}\bigr) \! + \!\frac{1}{2}{\delta^2}{\mathbb L}_0 f\bigl(t_{k-1}^{(l)},\mathbf{u}_{k-1}^{(l)}\bigr) \label{fun_fd}
\end{align}
and two differential operators ${\mathbb L}_0$ and ${\mathbb L}$ are defined by
\begin{align}
{\mathbb L}_0 & = \frac{\partial}{\partial t} + \sum \limits_{i = 1}^{K+1} f_i \frac{\partial}{\partial \mathbf{u}_i} + \frac{1}{2}\sum \limits_{j,p,r=1}^{K+1} G_{pj}  G_{rj} \frac{\partial^2}{\partial \mathbf{u}_p \partial \mathbf{u}_r}, \label{operator_L0} \\
{\mathbb L}_j &= \sum \limits_{i = 1}^{K+1} G_{ij} \frac{\partial}{\partial \mathbf{u}_i}, \quad i,j = 1, 2,\ldots, K+1. \label{operator_Lj}
\end{align}
where ${\mathbb L}$ stands for a square matrix with $(i,j)$ entry ${\mathbb L}_jf_i$.

Following~\cite{2022:Kulikova:ECC}, for the SDNF equation in~\eqref{eq1.3}, we have
\begin{align}
{\mathbb L}_0 f\bigl(t,\mathbf{u}\bigr) & = h_x\,\mathbf{V}\times \frac{\partial}{\partial t}I(\mathbf{x},t) + \sum \limits_{i = 1}^{K+1} f_i [J]_i
\end{align}
where $f_i$ is the $i$th entry of the vector-function in~\eqref{drift:state:fun} and $[J]_i$ is the  $i$th column of the matrix $J = \partial f\bigl(t,\mathbf{u}\bigr)/\partial \mathbf{u}$ in~\eqref{Jacobian:drift}.

The error covariance matrix is calculated as follows~\cite{2022:Kulikova:ECC}:
\begin{align}
& P_{k|k-1}^{(l+1)}  = \left[J_{fd}(t_k^{(l)},\mathbf{\hat u}_{k|k-1}^{(l)})\right] P_{k|k-1}^{(l)} \left[J_{fd}^{\top}(t_k^{(l)},\mathbf{\hat u}_{k|k-1}^{(l)})\right] \nonumber \\
& +\frac{\delta^2}{2}\left[G \;{\mathbb L} f^{\top}\!\!\bigl(t_{k-1}^{(l)},\mathbf{\hat u}_{k|k-1}^{(l)}\bigr) +{\mathbb L}  f\bigl(t_{k-1}^{(l)},\mathbf{\hat u}_{k|k-1}^{(l)}\bigr)\;G^{\top}\right] \nonumber \\
&+\!\!\frac{\delta^3}{3}\!\bigl[{\mathbb L} f\bigl(t_{k-1}^{(l)},\mathbf{\hat u}_{k|k-1}^{(l)}\bigr)\!\bigr]\bigl[{\mathbb L} f\bigl(t_{k-1}^{(l)},
\mathbf{\hat u}_{k|k-1}^{(l)}\bigr)\!\bigr]^{\top}\!\!\!+\!\delta G G^{\top}\label{EKF:covariance}
\end{align}
where the matrix ${\mathbb L}f\bigl(t,\mathbf{u}\bigr)$ is formed by columns as follows:
\begin{align}
{\mathbb L} f\bigl(t,\mathbf{u}\bigr) & = [\epsilon \mathbf{\Lambda}_{1,1}[J]_1 \: | \ldots | \: \epsilon \mathbf{\Lambda}_{(K+1),(K+1)}[J]_{K+1}]
\end{align}
 where $J_{fd} = \partial f\bigl(t,\mathbf{u}\bigr)/\partial \mathbf{u}$ is the Jacobian matrix derived for the discretized
drift function $f_d(\cdot)$ in~\eqref{fun_fd}. Next, the term $[J]_j$ stands for the $j$th column of the Jacobian matrix $J = \partial f\bigl(t,\mathbf{u}\bigr)/\partial \mathbf{u}$. Besides, the diffusion matrix is $G=\epsilon \mathbf{\Lambda}$, $\mathbf{\Lambda}:=\diag\{\lambda_0,\lambda_1,\ldots,\lambda_K\}$ with each $\lambda_k$, $k=0, \ldots, K$, computed by formula~\eqref{eq:lambdas}. The proof comes from formula~\eqref{operator_Lj} for calculating the operator ${\mathbb L}_j$. For the SDNF under examination, this  yields the columns ${\mathbb L}_j:=G_{j,j}[J]_j$. The readers are referred to~\cite{2022:Kulikova:ECC} for more details and derivation of the IT-1.5 EKF prediction step.

At the final iterate $L$ of the filter prediction step, set $\mathbf{\hat u}_{k|k-1}:=\mathbf{\hat u}^{(L)}_{k|k-1}$ and $P_{k|k-1} := P^{(L)}_{k|k-1}$.

\textsc{Measurement Update}: Apply formulas~\eqref{kf:f:K}~-- \eqref{kf:f:P} to get the filtered estimate $\mathbf{\hat u}_{k|k}$ and the error covariance $P_{k|k}$.

As shown in~\cite{2022:Kulikova:ECC}, the IT-1.5 EKF-based state reconstruction method outperforms the EM-0.5 EKF-based scheme for estimation accuracies and restoration quality of the SDNF membrane potential under examination. In the next section, we aim to investigate the pattern recognition facilities of two EKF-based techniques discussed.

\section{NUMERICAL EXPERIMENTS}\label{sect4}

To perform the comparative study of our two reconstruction methods under examination, we explore an example where a few active zones (they are called the bumps) appear in the neural tissue due to the external stimulus and the SDNF memory mechanism as discussed in~\cite{KuLi21CAM}. In this paper, we are interested in the quality of the pattern recognition in the neural tissue in case of the multi-bump solutions.

{\bf Example 1}. Let us consider the SDNF equation in~\eqref{eq1.3} with the decay rate $\alpha=1$, threshold $\theta = 0$ and the initial value $U(x,0)=U_0(x)=0$. The connectivity function is
\begin{align*}
F(|x-y|) & = A_1 e^{-k_1 |x-y|} \bigl[k_1 \sin(\pi |x-y|/k_2)  \\
& + \cos(\pi |x-y|/k_2)\bigr]
\end{align*}
where $A_1 = 2$, $k_1 = 0.08$ and $k_2 = 10$.

The problem is solved on the spatial domain $\Omega$ with the length $L = 100$, i.e. $x\in\Omega\equiv[-100, 100]$, and on the interval $t \in [0,10]$s. The external stimulus is given by
\begin{align}
I(x,t) & = -3.39967 + A e^{-(x-C)^2/2\sigma^2}, & t & \in [0, 5],  \label{ex2:stimulus:1}\\
I(x,t) & = -2.89967, & t & \in (5, 10] \label{ex2:stimulus:2}
\end{align}
where $A = 8$, $C = 0$. We examine $\sigma = 3$ and $\sigma = 13$ cases.

Following~\cite{KuLi21CAM}, Example~1 with $\sigma = 3$ yields an one-bump pattern in the neural tissue in case of a weak-noise case, say $\epsilon = 0.05$ in equation~\eqref{eq1.3}. Meanwhile, the external stimulus with $\sigma = 13$ and $\epsilon = 0.05$ provides a two-bump solution of the SDNF equation under examination. In our first set of numerical experiments, we illustrate the reconstruction process of the membrane potential in these two scenarios. For that, the following set of numerical experiments is performed. First, the SDNF equation is simulated with a small step size, $h_t=0.01$s, by using a high order method that is the the IT-1.5 integrator~\cite{1999:Kloeden:book}. We utilize $K=100$ eigenfunctions in the Karhunen-Lo\`eve formula~\eqref{eq2.8} where $\xi=0.1$ is used in~\eqref{eq:lambdas} for computing the diffusion matrix ${\mathbf \Lambda}$. The spatial discretization is performed with the step-size $h_x=0.2$ that yields $1000$ spatial points. The left graphs on Figs.~1 and~2 display the simulated reference or ``true'' solution of equation~\eqref{eq1.3} for the scenarios with $\sigma = 3$ and $\sigma = 13$, respectively. Second, given the ``true'' solution, the incomplete measurements are generated with the sampling rate $\Delta_t=0.5$s and spatial step size $\Delta_x=5$. The outcomes of this step are plotted on the intermediate graphs of Figs.~1 and~2, respectively. Finally, the membrane potential is recovered by the filtering method under discussion from the incomplete data available. For illustrative purpose, we utilize the EM-0.5 EKF reconstruction method from~\cite{2021:Kulikova:Romania} for the case of $\sigma = 3$ and plot the recovered membrane potential together with the data on the right plot of Fig.~1. Meanwhile, the scenario with $\sigma = 13$ is processed by the IT-1.5 EKF restoration process in~\cite{2022:Kulikova:ECC} and the results are illustrated by the right graph of Fig.~2. Both filtering methods are applied with $L=5$ subdivisions, the filter initial covariance $\Pi_0 = 0.1 I$ and the measurement noise covariance is assumed to be $R = 10^{-3} I$ where $I$ is an identity matrix of a proper size. Having analyzed the obtained results illustrated by Figs.~1 and~2, we observe a good reconstruction quality.

\begin{figure*}
\begin{tabular}{ m{0.3\textwidth}  m{0.3\textwidth}  m{0.3\textwidth} }
\includegraphics[width=0.3\textwidth]{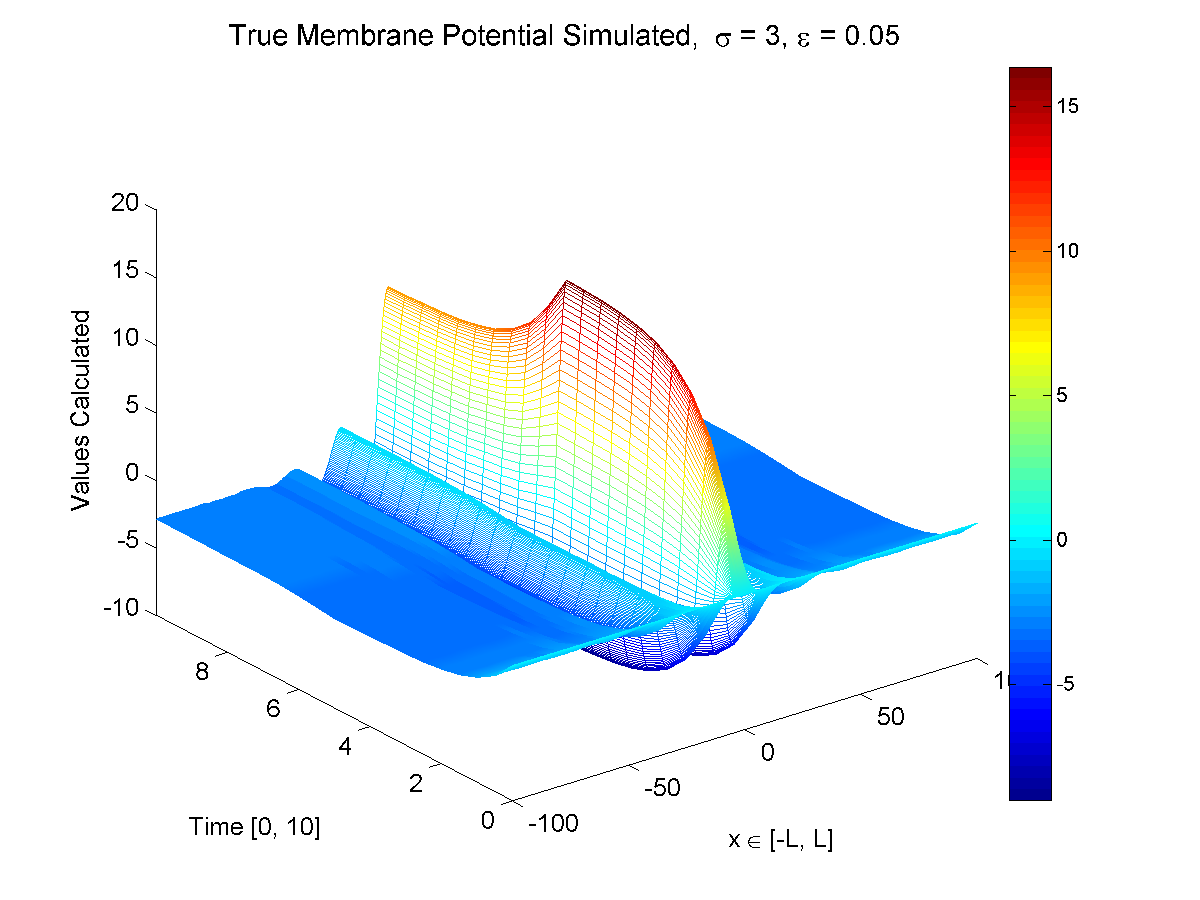} & \includegraphics[width=0.3\textwidth]{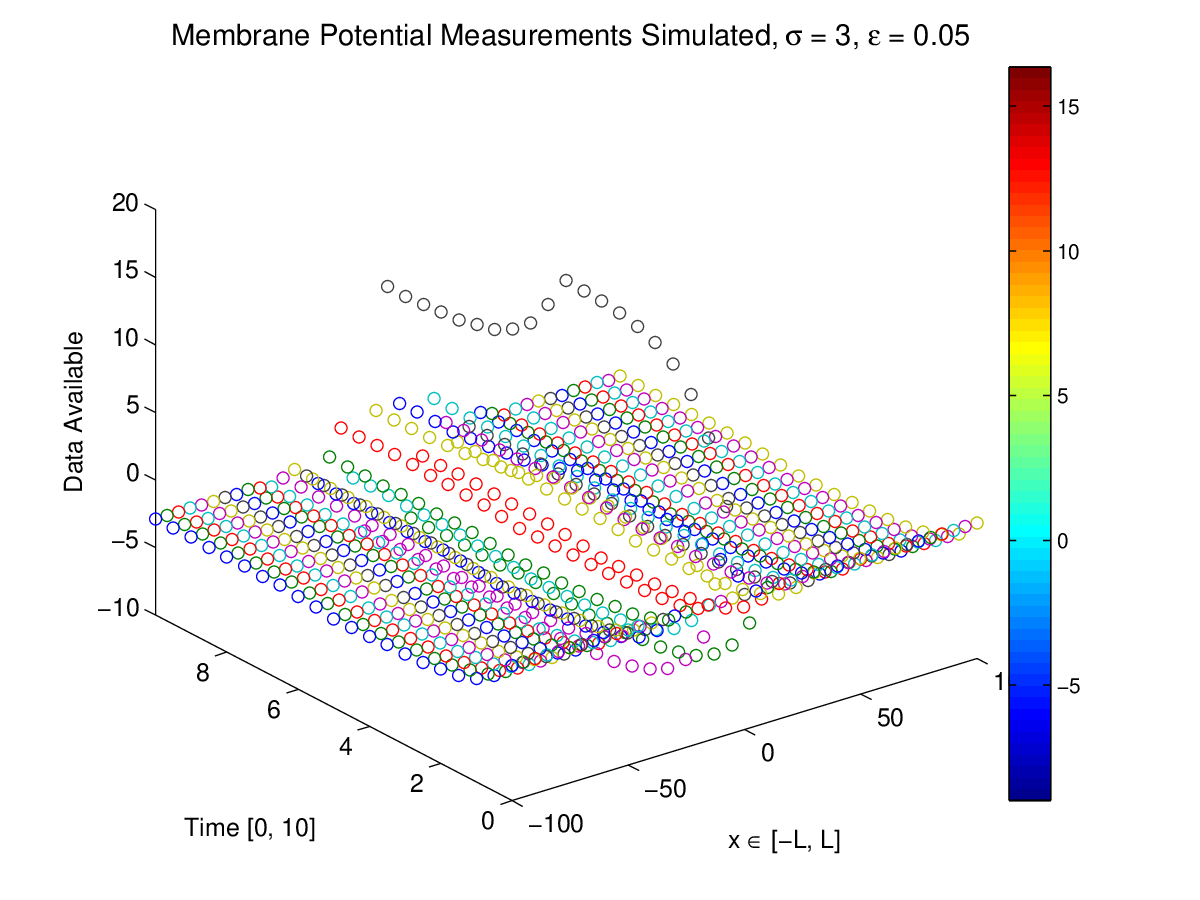} & \includegraphics[width=0.3\textwidth]{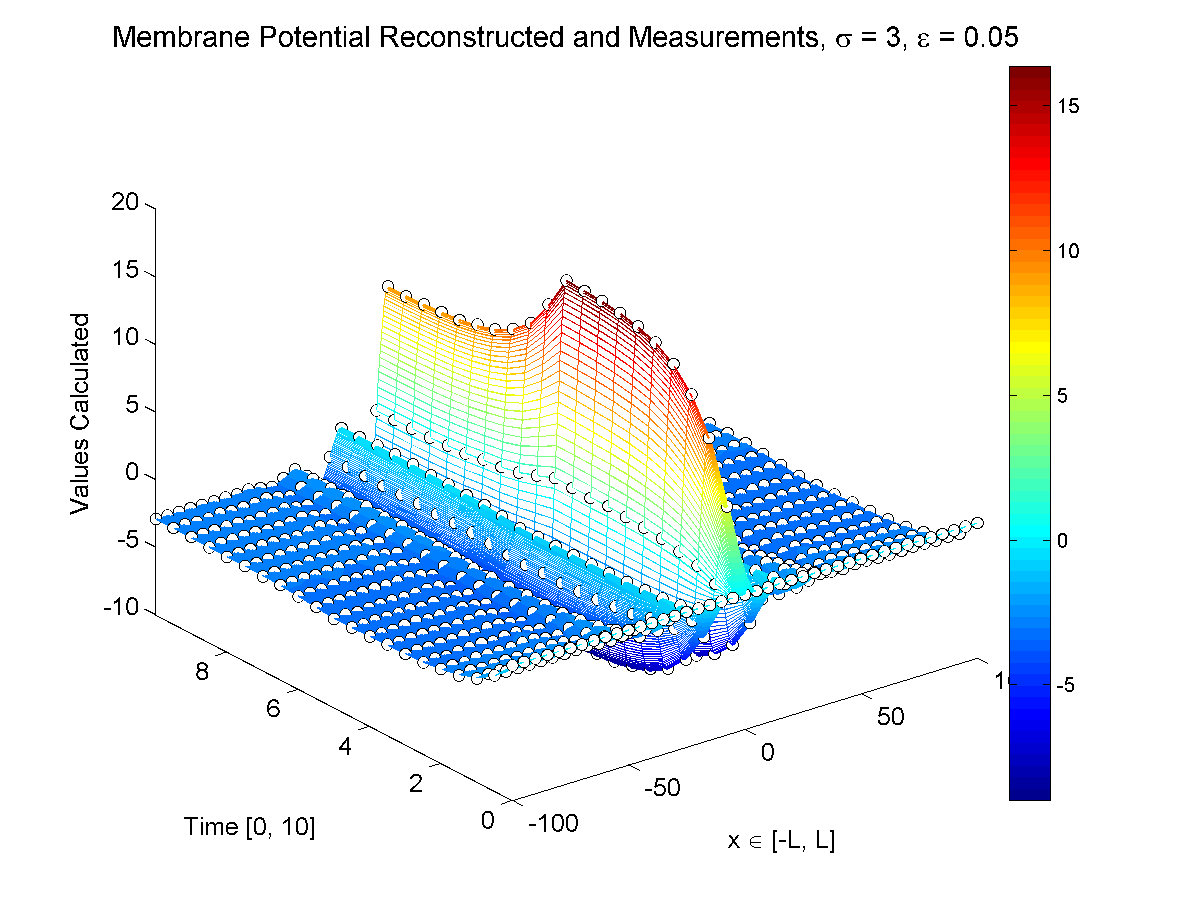} \\
\end{tabular}
\caption{The one-bump SNFE solution arisen in Example~1 when $\sigma = 3$ and the membrane potential reconstruction process.}
\end{figure*}

\begin{figure*}[th!]
\begin{tabular}{ m{0.3\textwidth} m{0.3\textwidth}  m{0.3\textwidth} }
\includegraphics[width=0.3\textwidth]{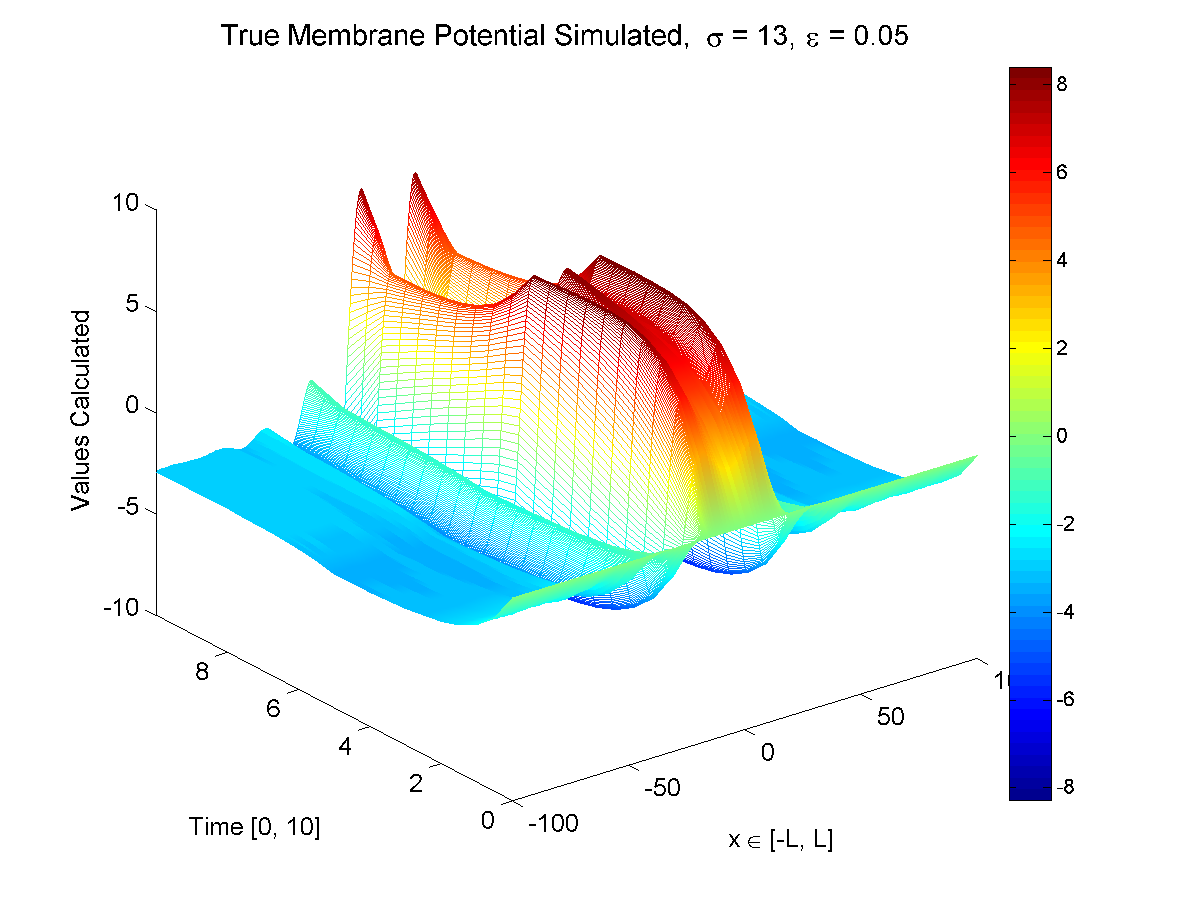} & \includegraphics[width=0.3\textwidth]{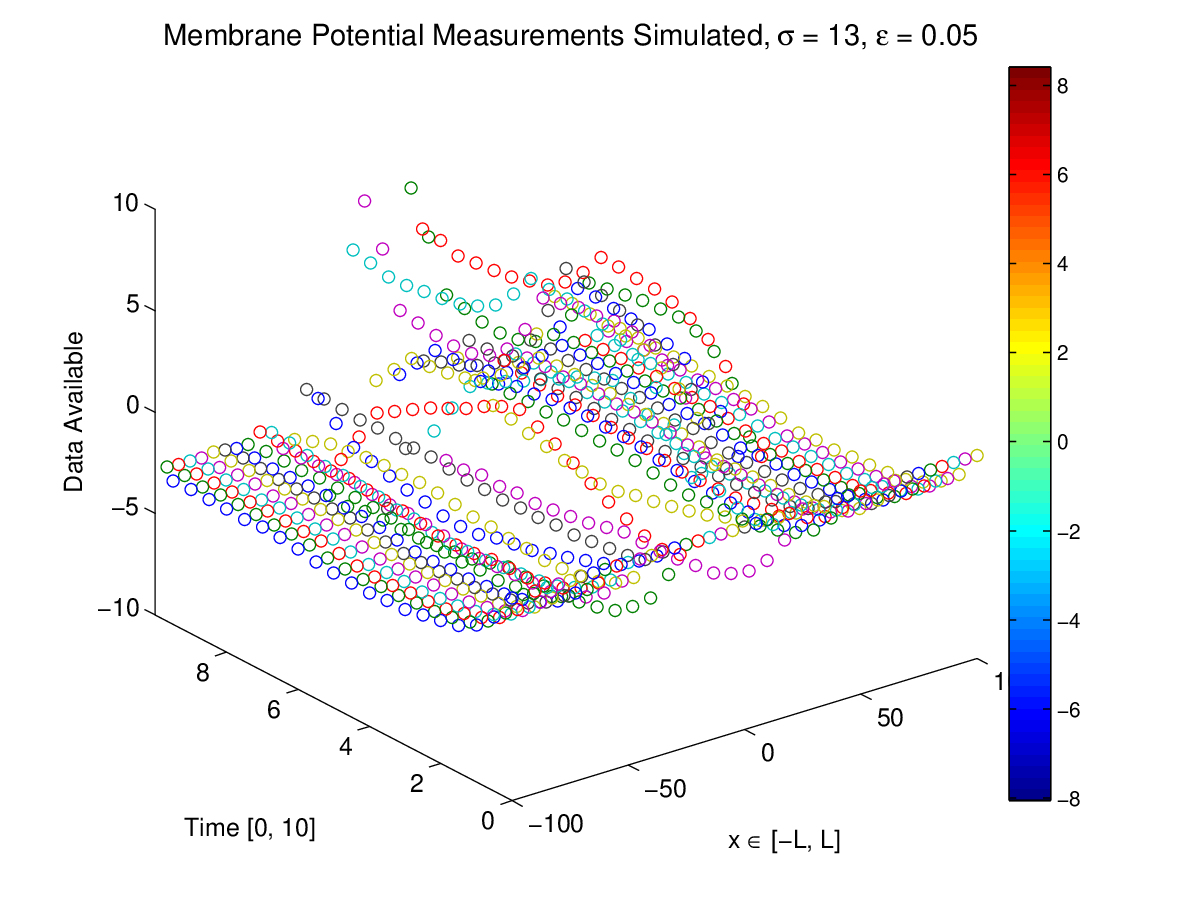} & \includegraphics[width=0.3\textwidth]{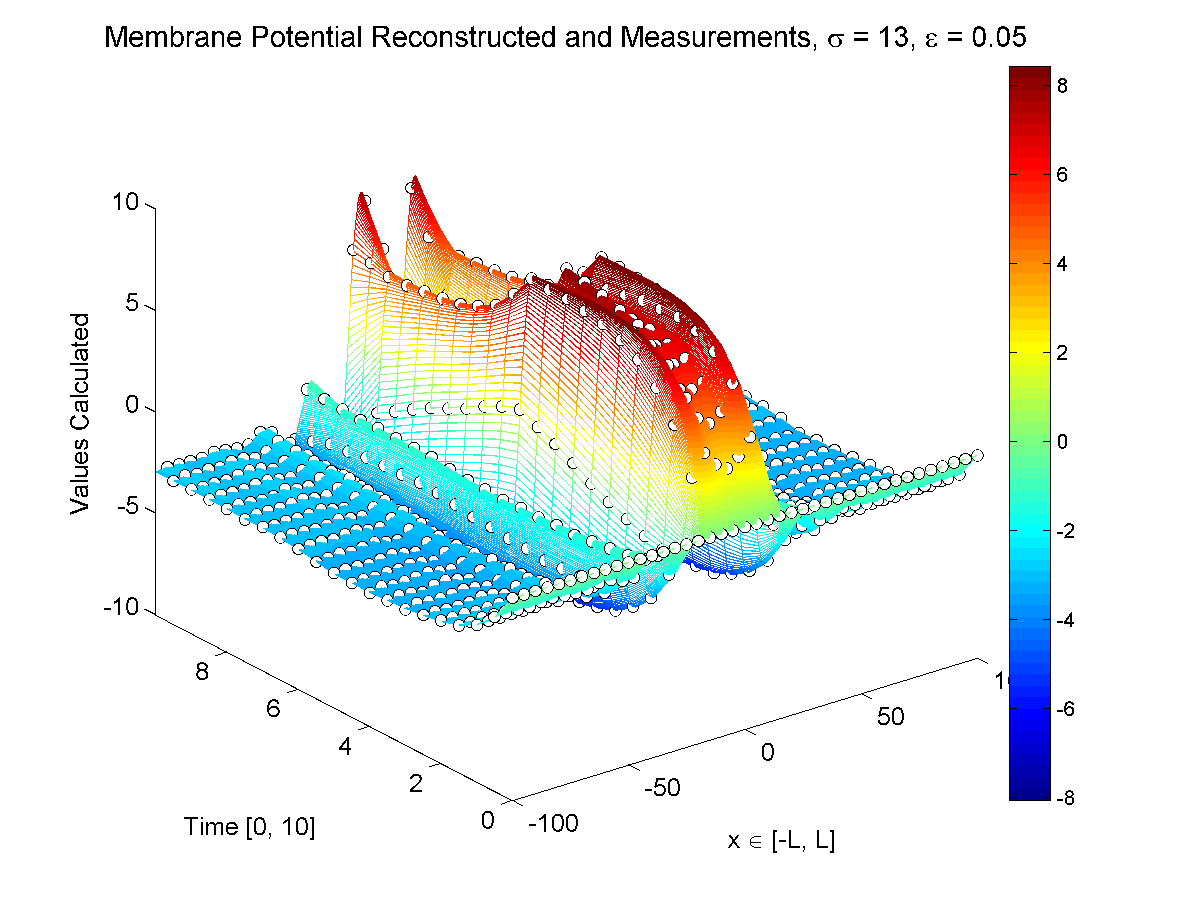} \\
\end{tabular}
\caption{The two-bump SNFE solution arisen in Example~1 when $\sigma = 13$ and the membrane potential reconstruction process.}
\end{figure*}

To get more insights about the restoration quality of two techniques under examination and to assess their reconstruction capacities, we next perform a set of numerical tests with $M=500$ Monte Carlo simulations as follows. First, the SDNF equation is simulated with $K=100$, $h_t=0.1$s, $h_x = 0.2$ for $\sigma = 3$ and $\sigma = 13$ and in case of two process noise levels $\epsilon = 0.05$ and $\epsilon = 0.5$. Second, the number of active zones (i.e. the bumps) at the last $T = 10$s are calculated out of $M=500$ trials. As discussed in~\cite{KuLi21CAM}, the multi-bump solutions are possible in the presence of model uncertainties, especially for the strong noise case scenario with $\epsilon = 0.5$. Third, the incomplete data sets are simulated in each Monte Carlo run with the sampling rate $\Delta_t = 0.2$s and the spatial step size $\Delta_x = 4$. This means that only one spatial point among $20$ is supplied by the sensor for collecting the data about the membrane potential. Besides, the measurements are taken at each $0.2$s. Next, given the simulated {\it incomplete} measurement history, we solve the inverse problem, i.e. the membrane potential is reconstructed from the incomplete data available. For that, we apply the EM-0.5 EKF method in~\cite{2021:Kulikova:Romania} and the IT-1.5 EKF scheme in~\cite{2022:Kulikova:ECC} with $L=1$ subdivisions and $\Pi_0 = 0.1 I$. Finally, we calculate the number of active zones recovered by the filtering methods at $T=10$s and compare these results with the number of bumps obtained from the ``true'' SDNF solution. This allows for calculating the number of mismatches in the pattern recognition process and, hence, to decide about the restoration quality and recognition facilities of two methods examined. The obtained results are collected for $\sigma = 3$ and $\sigma = 13$ cases in Tables~\ref{Tab:si3} and~\ref{Tab:si13}, respectively.

\begin{table*}
{\scriptsize
\caption{Pattern recognition quality of two EKF-based reconstruction methods in Example~1 with $\sigma = 3$ in $M=500$ trials.} \label{Tab:si3}
\centering
\begin{tabular}{c||c||c|c||c|c||c||c|c||c|c}
\hline
Number&  \multicolumn{5}{c||}{\bf Weak-noise case, $\epsilon =0.05$} &  \multicolumn{5}{c}{\bf Strong-noise case, $\epsilon =0.5$} \\
\cline{2-11}
of &  & \multicolumn{2}{c||}{\bf EM-0.5 EKF-based scheme} &
\multicolumn{2}{c||}{\bf IT-1.5 EKF-based scheme} & & \multicolumn{2}{c||}{\bf EM-0.5 EKF-based scheme} & \multicolumn{2}{c}{\bf IT-1.5 EKF-based scheme}\\
\cline{2-11}
Bumps &  Exact & Recovered & Mismatch &  Recovered & Mismatch & Exact & Recovered & Mismatch &  Recovered & Mismatch \\
\hline
\hline
1 & 500 &  500 & 0 & 500 & 0 & 458 & 457 & 1 & 457 & 1 \\
2 & 0   &  0   & 0 & 0   & 0 & 0   & 0   & 0 & 0   & 0 \\
3 & 0   &  0   & 0 & 0   & 0 & 38  & 39  & 1 & 39  & 1 \\
4 & 0   &  0   & 0 & 0   & 0 & 0   & 0   & 0 & 0   & 0 \\
5 & 0   &  0   & 0 & 0   & 0 & 4   & 4   & 0 & 4  & 0 \\
\hline
Total & 500   &  500   & 0 & 500   & 0 & 500  & 500  & 2 & 500  & 2 \\
\hline
\end{tabular}
}
\end{table*}

\begin{table*}
{\scriptsize
\caption{Pattern recognition quality of two EKF-based reconstruction methods in Example~1 with $\sigma = 13$ in $M=500$ trials.} \label{Tab:si13}
\centering
\begin{tabular}{c||c||c|c||c|c||c||c|c||c|c}
\hline
Number&  \multicolumn{5}{c||}{\bf Weak-noise case, $\epsilon =0.05$} &  \multicolumn{5}{c}{\bf Strong-noise case, $\epsilon =0.5$} \\
\cline{2-11}
of &  & \multicolumn{2}{c||}{\bf EM-0.5 EKF-based scheme} &
\multicolumn{2}{c||}{\bf IT-1.5 EKF-based scheme} & & \multicolumn{2}{c||}{\bf EM-0.5 EKF-based scheme} & \multicolumn{2}{c}{\bf IT-1.5 EKF-based scheme}\\
\cline{2-11}
Bumps &  Exact & Recovered & Mismatch &  Recovered & Mismatch & Exact & Recovered & Mismatch &  Recovered & Mismatch \\
\hline
\hline
1 & 0   &  0   & 0 & 0   & 0 & 60  & 55  & 5 & 55   & 5 \\
2 & 500 &  500 & 0 & 500 & 0 & 373 & 371 & 2 & 371  & 2 \\
3 & 0   &  0   & 0 & 0   & 0 & 48  & 52  & 4 & 52   & 4 \\
4 & 0   &  0   & 0 & 0   & 0 & 9   & 11  & 2 & 11   & 2 \\
5 & 0   &  0   & 0 & 0   & 0 & 10  & 11  & 1 & 11   & 1 \\
\hline
Total & 500   &  500   & 0 & 500   & 0 & 500  & 500  & 14 & 500  & 14 \\
\hline
\end{tabular}
}
\end{table*}

Having analyzed the results obtained, we observe a high pattern recognition quality of both filtering methods under examination. Indeed, both techniques carefully reconstruct the one-bump and two-bump membrane potential patterns (i.e. without mismatches in the number of bumps) in case of the weak-noise scenario with $\epsilon = 0.05$ for $\sigma = 3$ and $\sigma= 13$, respectively. As can be seen from the first panels of Tables~\ref{Tab:si3} and~\ref{Tab:si13} related to $\epsilon = 0.05$ case, there is no mismatch observed in the pattern recognition process in all $M=500$ trials. Next, the presence of the model uncertainties influences the pattern formation mechanism, significantly. This conclusion is inline with the results previously reported in~\cite{KuLi21CAM,Lima2022} and other papers. It is clearly seen form the second panel of Table~\ref{Tab:si3} that the one-bump ``true'' solution might be destroyed because of the presence of model uncertainties with the noise level $\epsilon = 0.5$. This yields the multi-bump patterns in the neural tissue. A similar conclusion is made for the two-bump solution related to the scenario with $\sigma = 13$. We also conclude that both reconstruction methods recover the pattern formation process in case of a strong noise scenario with a good recognition quality. Indeed, the total number of mismatches does not exceed 3\% of the total $M=500$ variants. As mentioned above, in case of the weak-noise scenario they recover the ``exact'' neural tissue pattern without mismatches at all.

Although the IT-1.5 EKF-based reconstruction method is shown to be more accurate than the EM-0.5 EKF-based scheme~\cite{2022:Kulikova:ECC}, we observe their  identical performance  in terms of the recognition facilities. Following the cited paper, a significant difference in their reconstruction quality is observed when the spatial step size $\Delta_x$ increases, i.e. when less measurements are available on the domain $[-L, L]$. Therefore, to scrutinize their performance, we repeat the numerical experiments described above for various values of $\Delta_x$  and calculate the total number of mismatches in $M=1000$ simulations. We again observed their equal performance with very high reconstruction quality in case of the weak-noise scenario with $\epsilon = 0.05$, i.e. both methods recognize the patterns in the neural tissue without any mismatch. Meanwhile, the strong-noise case $\epsilon = 0.5$ yields the results illustrated by Figs.~3 and~4. It is clearly seen that the IT-1.5 EKF-based method outperforms the EM-0.5 EKF alternative for its recognition facilities because it yields a less number of mismatches in case of decreasing spatial information available about the membrane potential.

\begin{figure}[th!]
\includegraphics[width=0.46\textwidth]{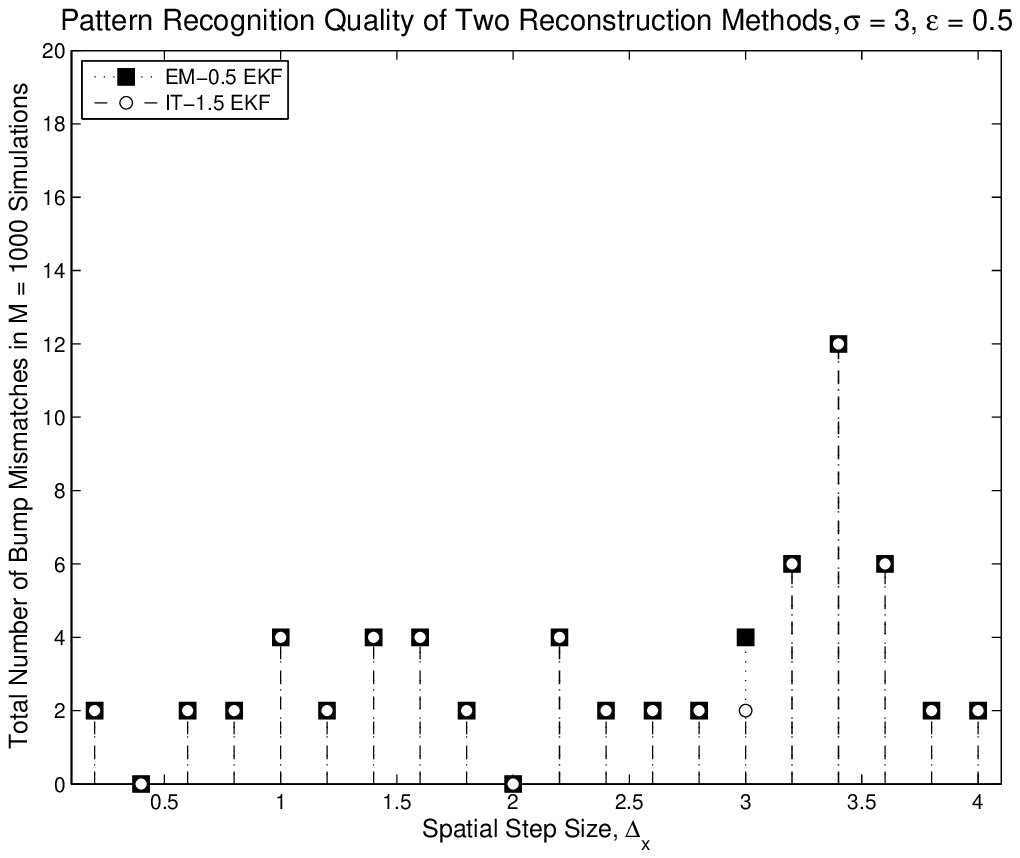}
\caption{Recognition facilities of two EKF-based reconstruction methods. }
\end{figure}

\begin{figure}[th!]
\includegraphics[width=0.46\textwidth]{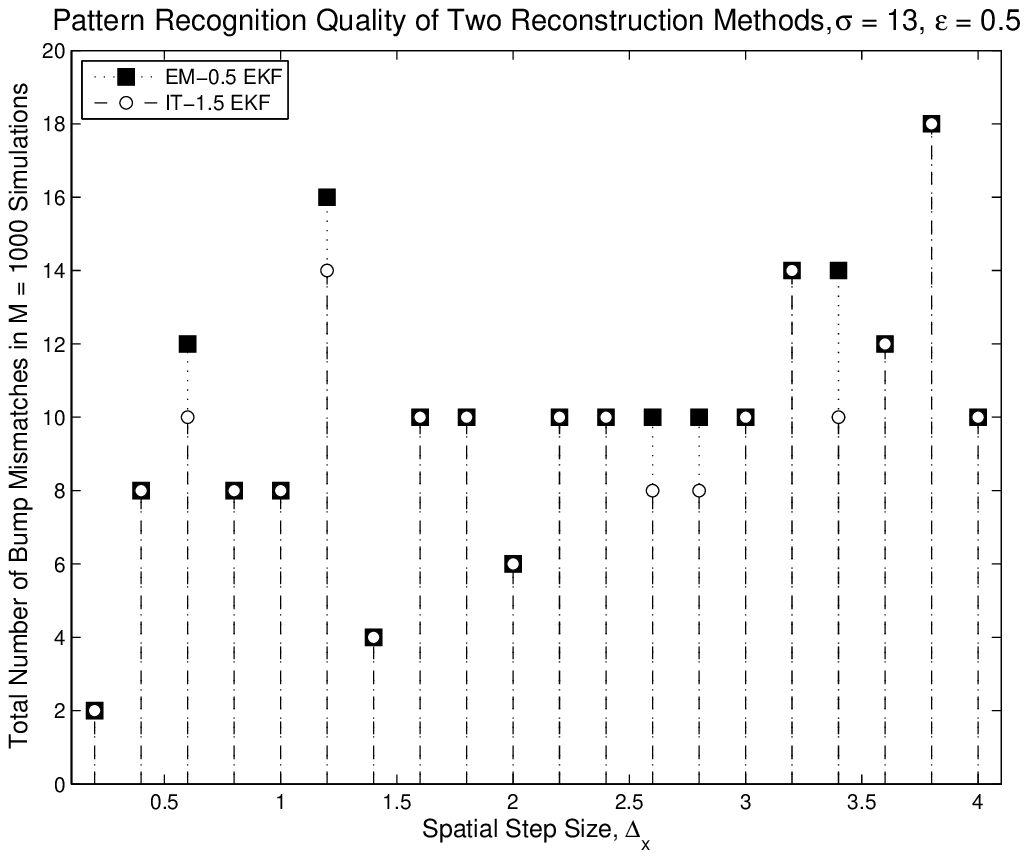}
\caption{Recognition facilities of two EKF-based reconstruction methods. }
\end{figure}

\bibliographystyle{IEEEtran}
\bibliography{IEEEabrv,Library/books,%
              Library/list_kf_pde,%
              Library/list_nonlinear,%
              Library/list_Kulikov,%
              Library/list_neural}

\end{document}